\documentstyle[12pt,dina4,fleqn,emlines,amsfonts]{article}
\newcommand{\be}{\begin{enumerate}}
\newcommand{\ol}[1]{\overline{#1}}

\newcommand{\ee}{\end{enumerate}}
\newcommand{\bq}{\begin{quote}}
\newcommand{\eq}{\end{quote}}
\newcommand{\bi}{\begin{itemize}}
\newcommand{\ei}{\end{itemize}}

\begin{document}
\begin{center}
{\Large{\bf On a strange observation in the theory of the dimer problem}}\\[4em]
Peter E. John and Horst Sachs\\
Technical University of Ilmenau,\\
Institute of Mathematics,\\
Ilmenau, Germany\\
\end{center}
\vspace{0,5cm}
Report-no: TUI-Math 16/97\\
Subj-class: Combinatorics, Graph Theory\\
MSC-class: 05B45 (Primary) 05A15, 05C70 (Secondary)\\[2em]

This is a contribution to the number theory of the dimer problem.

Let $A_n = 2^n \cdot B_n^2~ (B_n > 0)$ denote the number of dimer coverings 
(i.e., perfect matchings) of a $2n \times 2n$ square lattice graph. $B_n$ 
turns out to be
an integer. Motivated by a somewhat strange observation (see below), we
investigated the residue classes $B_n~ mod~ 2^r$. In this paper, we outline
a method that, for a fixed integer $r$, enables these residue classes to be
determined. Explicitly, $B_n~ mod~ 64$ is calculated and from the resulting
formula the following proposition is deduced:
$$B_n \equiv \left\{ \begin{array}{lll}
n+1 & (\mbox{mod}~ 32) & if~ n~ is~ even\\
(-1)^{\frac{n-1}{2}} \cdot n & (\mbox{mod}~ 32) & if~ n~ is~ odd.
\end{array} \right.$$
The analogous statement $mod~64$ is not true: $B_3 = 29$ is a counterexample.\\[1em] 
{\bf 1. Introduction.}\\[1em] Let $A_n$ denote the number of different ways in
which a $2n \times 2n$ chessboard can be covered with $2n^2$ dominoes such
that each domino covers two adjacent squares $(n = 1,2,3,...)$. Set\\[1em]
(1)~~~~~~~~~~~ $A_n = 2^nB_n^2$~~~~~where~~ $B_n > 0$.\\[1em]
As we shall see (and is well known), $B_n$ is an integer.\\ Let
$$u_k = 2 \mbox{cos} \frac{k \pi}{2n+1}~,~~ k = 0, \pm 1, \pm 2, ...;$$
note that\\[1em]
(2)~~~~~~~~~~~ $u_{-k} = u_k,~~ u_{2n+1-k} = - u_k$.\footnotemark
\footnotetext{In fact, $u_k$ depends on two variables, $n$ and $k$, but for the sake of
readability we shall suppress the $n$.
The same applies to the following
function symbols used in this paper: $v_k, \ol{v}_k, s_k, \ol{s}_k, \sigma_k.$}\\

The following is a classical result of Fisher/Temperley and Kasteleyn, 1961
(see, e.g. [2], p. 96):
\begin{center}
$A_n = \prod \limits_{j=1}^n \prod \limits_{k=1}^n (u_j^2 + u_k^2).$\\[1em]
\end{center}
Some first investigations of the numbers $B_n$ (for small value of $n$)
 gave rise to a somewhat strange
conjecture which was backed by the list of the first 20 numbers $B_n$ provided
by N. Saldanha (Lyon)\footnotemark 
\footnotetext{e-mail message of May 9, 1997}; this conjecture,
now Theorem A, is to be proved in this paper:\\[1em]
{\bf Theorem A.} $The~ statement~ {\Bbb S}_r:$
$$B_n \equiv \left\{ \begin{array}{rl}
n+1 & if~ n~ is~ even\\
(-1)^{\frac{n-1}{2}} \cdot n & if~ n~ is~ odd
\end{array} \right. (\mbox{mod}~2^r)$$
$is~ true~ for~ r = 5~ (thus~ also~ for~ r = 1,2,3,4)~ and~ false~ for~ r = 6~
(thus~ also~ for\\ r > 6).$\\[1em]
The numbers $B_3 = 29,~ B_5 = 8	9893
\equiv 37$ and $B_6 = 28793575 \equiv 39,~ mod~64$ are counterexamples
to statement ${\Bbb S}_6$.\\ \\ 
An immediate consequence of Theorem A is\\[1em]
{\bf Corollary A.} ([1]). $All~ numbers~ B_n~ are~ odd.$\\

We shall prove a slightly more general proposition, namely, we shall 
establish a formula for the residue classes $B_n~mod~64$ (Theorem B). The 
method used can, in principle, be extended so as to enable the 
$B_n~mod~2^r$ to be
calculated for any fixed positive integer $r$, however, with increasing $r$
the situation becomes more and more involved: we are unable to provide a
closed formula (in $n$ and $r$) for $B_n~mod~ 2^r$.\\[1em]
{\bf 2. Preliminaries.}\\[1em]
Let\\ \\
(3)~~~~~~~~~~~ $v_k = u_{2k} = 2~\mbox{cos}~\frac{2k\pi}{2n+1}$;\\ \\
note that\\ \\
(4)~~~~~~~~~~~ $v_{-k} = v_k = v_{2n+1-k}$.\\ \\
By (2), $u_{2k+1} = - u_{2(n-k)},$
thus, by (3),
$$u_{2k} = v_k,~ u_{2k+1} = - v_{n-k}.$$
This implies
$$\{u_k^2~ |~ k = 1,2,...,n\} = \{v_k^2~ |~ k = 1,2,...,n\},$$
therefore,\\ \\
(5)~~~~~~~~~~~ $A_n = \prod \limits_{j=1}^{n} \prod \limits_{k=1}^{n}
   (v_j^2 + v_k^2).$\\ \\
Define
$$\begin{array}{rcl}
2 \alpha & = & 2 \alpha(x)~ =~ \mbox{arc~ cos}~ \frac{x}{2},\\ \\
S_n = S_n(x) & = & \displaystyle \frac{\mbox{sin}~ [(2n+1) \alpha(x)]}{\mbox{sin}~\alpha(x)}~~~~~~~~~~ (-2 < x < 2).
\end{array}$$
Note that\\ \\
(6)~~~~~~~ $2~ \mbox{cos}(2\alpha(x)) = x,~~\mbox{sin}~ \alpha(x)
    = \frac{1}{2} \sqrt{2-x},~~ S_{-1}(x) = -1,~~ S_0(x) = 1.$\\[1em]
Using (6), from\\ \\
$\mbox{sin}~ \alpha \cdot S_{n \pm 1} = \mbox{sin} [(2n + 1 \pm 2)\alpha] =
 \mbox{cos}~(2 \alpha) \cdot \mbox{sin} [(2n+1) \alpha] \pm \mbox{sin} (2 \alpha) 
 \cdot \mbox{cos}[(2n+1) \alpha]$\\ \\
we obtain
$$S_{n+1} + S_{n-1} = 2~ \mbox{cos}(2 \alpha) \frac{\mbox{sin}[(2n+1) 
 \alpha]}{\mbox{sin}~ \alpha}
= x \cdot S_n .$$
Thus $S_n$ satisfies the recurrence relation 
$$S_{n+1} = x~S_n - S_{n-1};~ S_{-1} = -1,~ S_0 = 1.$$
This implies that $S_n(x)~~ (n = 0,1,2,...)$ is a polynomial in $x$ of degree
$n$ with integral rational coefficients and first term $x^n$. Set
\begin{center}
$S_n(x) = \sum \limits_{\nu = 0}^{n} (-1)^{\nu} s_{\nu} x^{n-\nu},~~~
n = 0,1,2,...~~ .$
\end{center}
It is easily verified by induction that\\[1em]
$(7)~~~~~~~ \left \{ \begin{array}{lcl}
s_{2 \rho}   & = & (-1)^{\rho} {n-\rho \choose \rho},\\ \\
s_{2 \rho+1} & = & (-1)^{\rho+1} {n-1-\rho \choose \rho}.
\end{array} \right.$\\[1em]
Further, $\alpha(v_k) = \frac{k \pi}{2n+1}$, thus
$$S_n(v_k) = \frac{1}{\sqrt{2-v_k}} \mbox{sin}(k \pi) = 0,~~~~~ k = 1,2,...,n.$$
This implies that
\begin{center}
$S_n(x) = \prod \limits_{k=1}^{n} (x - v_k) = \sum \limits_{\nu = 0}^{n}
(-1)^{\nu} s_{\nu}~ x^{n-\nu}.$
\end{center}
Note that, by (7),\\[1em]
$(8)~~~~~~~ \prod \limits_{k=1}^{n} v_k = s_n = \left\{ \begin{array}{rl}
(-1)^{\frac{n}{2}} & \mbox{if $n$ is even,}\\
(-1)^{\frac{n+1}{2}} & \mbox{if $n$ is odd.}
\end{array} \right.$\\[1em]
Equation (5) implies\\[1em]
(9)~~~~~~~ $A_n = \prod \limits_{i=1}^{n} (2v_i^2)~ \{ \prod \limits_{(j,k)}
       (v_j^2 + v_k^2)\}^2$\\[1em]
where the symbol $(j,k)$ indicates that the operation to which it pertains is
performed over the set ${\Bbb P}$ of all ${n \choose 2}$ unordered pairs 
$(j,k),~~ j,k \in \{1,2,...,n\}, ~ j \not= k$.\\

By (8), $\prod \limits_{i=1}^n v_i^2 = s_n^2 = 1$, thus, comparing (1) and
(9), we have obtained\\[1em]
(10)~~~~~~~ $B_n = \prod \limits_{(j,k)} (v_j^2 + v_k^2)$.\\[1em]
As the last product is an integral symmetric function of the $v_k$, $B_n$
is rational as well as integral.\\[1em]
{\bf 3. Proof of Theorem A.}\\[1em]
From the trigonometric identity
$$\mbox{cos}^2 x + \mbox{cos}^2 y = 1 + \mbox{cos} (x + y)~ \mbox{cos} (x-y)$$
we obtain\\[1em]
(11)~~~~~~~ $v_j^2 + v_k^2 = 4 + v_{k+j} v_{k-j}.$\\[1em]
Let $w_k$ denote an indeterminate subject to the rule\\[1em]
(12)~~~~~~~ $w_{-k} = w_k = w_{2n+1-k}~~ (k = 1,2,...,n);$\\[1em]
note that, by (4), the $v_k$ satisfy (12).\\

It is not difficult to prove that\\[1em]
(13)~~~~~~~ $\{w_{k+j} w_{k-j}~ |~ (j,k) \in {\Bbb P}\} = \{ w_j w_k~ |~ (j,k) \in
{\Bbb P}\}.$\\[1em]
Using (11), define\\[1em]
$\begin{array}{lcl}
(14)~~~~~~~ h_n(x) & = & \prod \limits_{(j,k)} (x + v_{k+j} v_{k-j}) 
               =  \prod \limits_{(j,k)} (x - 4 + v_j^2 + v_k^2)\\
              & = & x^q + \sigma_1 x^{q-1} + ... + \sigma_q~~~
                   \mbox{where}~ q = {n \choose 2}.
\end{array}$\\[1em]
Clearly, the $\sigma_k$ are integral rationals.\\

Set\\[1em]  
(15)~~~~~~~ $U_n = \sigma_q,~ V_n = U_n G_n = \sigma_{q-1},~
W_n = U_n H_n = \sigma_{q-2}.$\\[1em]
By (10), (14) and (15),\\[1em]
(16)~~~~~~~ $B_n = h_n(4) \equiv 16 W_n + 4 V_n + U_n
                \equiv  U_n (16 H_n + 4G_n + 1) ~~ (\mbox{mod}~64).$\\[1em]
We shall in order calculate $U_n,~G_n,~H_n$.
\begin{center}
$\star$
\end{center}
By (14) and (15),\\[1em]
(17)~~~~~~~ $U_n = \sigma_q = \prod \limits_{(j,k)} (v_{k+j}~ v_{k-j})$.\\[1em]
Applying (13) (with $w_k = v_k$), we have\\[1em]
$(18)~~~~~~~ U_n = \prod \limits_{(j,k)} (v_j v_k) = \{\prod \limits_{i=1}^n
 v_i\}^{n-1} = \left\{ \begin{array}{cr}
(-1)^{\frac{n}{2}} & \mbox~{if~ n~ is~ even}\\
     1             & \mbox~{if~ n~ is~ odd.}
\end{array} \right.$\\[1em]
Note that this formula implies the validity of statements ${\Bbb S}_1$ and
${\Bbb S}_2$.
\begin{center}
$\star$
\end{center}  
Set $\ol{v}_k = v^{-1}_k$. Clearly, by (4), the $\ol{v}_k$ satisfy (12).
By (15) and (14),\\
\begin{center}
$U_n G_n = V_n = \sigma_{q-1} = \sum \limits_{(j,k)} \frac{\sigma_q}
{v_{k+j} v_{k-j}} = U_n \sum \limits_{(j,k)} \ol{v}_{k+j} \ol{v}_{k-j},$
\end{center}
thus, using (13) (with $w_k = \ol{v}_k$),\\[1em]
$(19)~~~~~~~ G_n = \sum \limits_{(j,k)} \ol{v}_{k+j}~ \ol{v}_{k-j} =
  \sum \limits_{(j,k)} \ol{v}_j~\ol{v}_k.$\\[1em]
The numbers $\ol{v}_k$ are the zeros of the ``reverse'' polynomial of
$S_n(x)$, namely,\\[1em]
$(20)~~~~~~~ \ol{S}_n(x) = (-1)^n s^{-1}_n x^n S_n (\frac{1}{x}) = x^n - \ol{s}_1
x^{n-1} + - ... + (-1)^n \ol{s}_n$\\[1em]
where\\[1em]
$(21)~~~~~~~ \ol{s}_k = s^{-1}_n s_{n-k},~~ k = 1,2,...,n.$\\[1em]
By (19), (20) and (21),\\
$$G_n = \ol{s}_2 = s_n^{-1} s_{n-2};$$
therefore, by virtue of (7),\\[1em]
$(22)~~~~~~~ G_{2 \lambda} = G_{2 \lambda+1} = - {\lambda +1 \choose 2}.$\\[1em]
Formul{\ae} (16), (18) and (22) imply the validity of statements ${\Bbb S}_3$
and ${\Bbb S}_4$.
\begin{center}
$\star$
\end{center}
By (15) and (17),\\

$U_n H_n = W_n = \sigma_{q-2} = \sum \limits_{((i,j),(k,l))}
 \frac{\sigma_q}{v_{j+i} v_{j-i} v_{l+k} v_{l-k}}
         = U_n \sum \limits_{((i,j),(k,l))} \ol{v}_{j+i} \ol{v}_{j-i}
              \ol{v}_{l+k} \ol{v}_{l-k}.$\\[1em]
(13) (with $w_k = \ol{v}_k$) now yields
\begin{center}
$H_n = \sum \limits_{((i,j),(k,l))} \ol{v}_i \ol{v}_j \ol{v}_k \ol{v}_l.$
\end{center}
The last sum consists of ${{n \choose 2} \choose 2}$ terms, each of the form
$\ol{v}_p \ol{v}_q \ol{v}_r \ol{v}_s$ (type 1) or 
$\ol{v}_p^2 \ol{v}_q \ol{v}_r$ (type 2) where the $p,q,r,s$ are pairwise
distinct. Each term of type 1 occurs precisely thrice since a set of 4 
elements can be partitioned into two pairs in exactly 3 ways; each term of
type 2 occurs precisely once since a family $\{a,a,b,c\}$ $(a \not= b \not= c
\not= a)$ can be partitioned into two pairs $(x,y)$, $x \not= y$, in exactly
one way. Thus the number of terms of type 1 and of type 2 is $3 {n \choose 4}$
and $3 {n \choose 3}$, respectively, which add up to the total number of 
terms, namely, $3 {n+1 \choose 4} = {{n \choose 2} \choose 2}$.\\

Denote the symmetric power sum in the variables $x_1, x_2,..., x_n$ that is
the sum of all terms derivable from $x_1^{e_1} x_2^{e_2} ... x_t^{e_t}$
$(e_{\tau} > 0)$ by permuting the $x_1, x_2, ..., x_n$ (without repetition)
by $[x_1^{e_1} x_2^{e_2} ... x_t^{e_t}]_n$ \footnote {E.g., $[x_1]_4 = x_1 + x_2
+ x_3 + x_4$, $[x_1^2 x_2]_3 = x_1^2 x_2 + x_1^2 x_3 + x_2^2x_1 + x_2^2 x_3 +
x_3^2 x_1 + x_3^2 x_2$.}. In this notation,\\[1em]
$(23)~~~~~~~ H_n = 3[\ol{v}_1 \ol{v}_2 \ol{v}_3 \ol{v}_4]_n + [\ol{v}_1^2 \ol{v}_2
\ol{v}_3]_n$.\\[1em]
It can easily be checked that\\[1em]
$(24)~~~~~~~ [x_1]_n~ [x_1 x_2 x_3]_n = [x_1^2 x_2 x_3]_n + 4 [x_1 x_2 x_3
x_4]_n,$\\[1em]
thus, by (23), (24) and (21),\\
\begin{eqnarray*}
H_n & = & [\ol{v}_1]_n [\ol{v}_1 \ol{v}_2 \ol{v}_3]_n -  
[\ol{v}_1 \ol{v}_2 \ol{v}_3 \ol{v}_4]_n\\
    & = & \ol{s}~ \ol{s}_3 - \ol{s}_4 = s_n^{-2} (s_{n-1} s_{n-3} - s_n s_{n-4})
    = s_{n-1} s_{n-3} - s_n s_{n-4}.
\end{eqnarray*}\\
(7) now yields\\[1em]
$(25)~~~~~~~ \left\{ \begin{array}{lcl}
H_{2 \lambda} & = & - \{ \lambda {\lambda + 1 \choose 3} + {\lambda + 2 
\choose 4}\},\\
H_{2 \lambda + 1} & = &  - \{ (\lambda + 1) {\lambda + 2 \choose 3} + {\lambda + 2 
\choose 4}\}.
\end{array} \right.$\\[1em]
\begin{center}
$\star$
\end{center}
(16), (18), (22) and (25) imply\\ \\
{\bf Theorem B.}\\
$$\begin{array}{lcl}
B_{2 \lambda + 1} & \equiv & - 16 (\lambda + 1) {\lambda + 2 \choose 3} - 16 {\lambda + 2 
\choose 4} - 4 {\lambda + 1 \choose 2} + 1~~~~ (\mbox{mod}~ 64),\\ \\
(-1)^{\lambda} B_{2 \lambda} & \equiv &  - 16 \lambda {\lambda + 1 \choose 3}
- 16 {\lambda + 2 \choose 4} - 4 {\lambda + 1 \choose 2} + 1~~~~ (\mbox{mod}~ 64).
\end{array}$$
\begin{flushright}
$\Box$
\end{flushright}
Replacing $\lambda$ by $\lambda + 16$, we obtain\\[1em]
{\bf Corollary 1.}~~~ $B_{n+32} \equiv B_n + 32~~~~ (\mbox{mod}~ 64) \hfill \Box$\\[1em]
implying\\[1em]
{\bf Corollary 2.}~~~ $B_n~~ \mbox{mod}~ 64~ is~ periodic~ in~ n~ with~
 primitive~ period~ 64. \hfill \Box$\\

Note that in the congruences of Theorem B the first terms on the right-hand side are
divisible by 32, thus\\

$B_{2 \lambda + 1} \equiv (-1)^{\lambda} B_{2 \lambda} \equiv
 - 16 {\lambda + 2 \choose 4} - 4 {\lambda + 1 \choose 2} + 1~~~~ 
(\mbox{mod}~ 32).$\\[1em]
Set~ $\lambda = 8 \kappa + \vartheta,~~~~~ \vartheta = 0,1,2,...,7$.\\
Note that ${\lambda + 2 \choose 4}$ is even if $\vartheta \in \{0,1,6,7\}$
and odd otherwise. This implies\\

$B_{2\lambda+1} \equiv (-1)^{\vartheta} B_{2\lambda} \equiv
\left\{ \begin{array}{lcl}
- 2 \lambda (\lambda + 1) + 1~~ (\mbox{mod}~ 32) & if & \vartheta \in \{0,1,6,7\}\\
- 2 \lambda (\lambda + 1) + 17 ~ (\mbox{mod}~ 32) & if & \vartheta \in \{2,3,4,5\}. 
\end{array} \right.$\\[1em]
After an easy calculation we obtain\\

$B_{2\lambda+1} \equiv (-1)^{\vartheta} B_{2\lambda} \equiv (-1)^{\vartheta}
(2\lambda + 1)~~~~ (\mbox{mod}~ 32).$\\[1em]
This implies the validity of statement ${\Bbb S}_5$ and thus proves
Theorem A. \hfill $\Box$\\ \\
{\bf 4. Concluding Remark.}\\[1em]
Note that the numbers $v_k^2$ are pairwise distinct. Set
$P_n(x) = \prod \limits_{j=1}^n (x - v_j^2)$. Clearly,
$P_n(x^2) = \prod \limits_{j=1}^n (x - v_j) \prod \limits_{k=1}^n (x + v_k) =
(-1)^n S_n (x) S_n (-x)$, therefore,\\ 
\begin{eqnarray*}
\frac{1}{4} \sqrt{4-x^2} P_n(x^2) & = & (-1)^n \sin[(2n+1) \alpha(x)] \cdot
\sin[(2n+1) \alpha(-x)]\\
& = & \sin[(2n+1) \alpha(x)] \cdot \cos [(2n+1) \alpha (x)]\\
& = & \frac{1}{2} \sin [(2n+1)~\mbox{arc}~ \cos \frac{x}{2}]
 =  \frac{1}{4} \sqrt{4-x^2}~ {\cal U}_{2n} (\frac{x}{2})
\end{eqnarray*}
where\\

${\cal U}_{2n}(x) = \sum \limits_{\nu = 0}^n (-1)^{\nu} {2n-\nu \choose \nu}
(2x)^{2(n-\nu)}$\\[1em]
is the (2$n$)th Chebyshev polynomial of the second kind. Thus\\

$P_n(x) = \sum \limits_{\nu = 0}^n (-1)^{\nu} {2n-\nu \choose \nu} x^{n-\nu}.$\\[1em]
The numbers $v_k^4$ are the zeros of the polynomial $Q_n(x)$ defined by\\

$Q_n(x^2) = (-1)^n P_n (x) P_n (-x).$\\[1em]
The discriminants of $P_n$ and $Q_n$ are\\

$\Delta(P_n) = \prod \limits_{(j,k)} (v_j^2 - v_k^2)^2$~~ and~~
$\Delta(Q_n) = \prod \limits_{(j,k)} (v_j^4 - v_k^4)^2$,\\[1em]
respectively, and we obtain\\

$B_n^2 = \prod \limits_{(j,k)} (v_j^2 + v_k^2)^2 = \displaystyle \frac{\Delta(Q_n)}
{\Delta(P_n)}.$\\[1em]
Thus the problem of expressing $B_n$ through the coefficients of $S_n$, or
of $P_n$, is closely related to the task of expressing the discriminant of
a polynomial $F$ through the coefficients of $F$.\\ \\
{\bf References.}
\begin{itemize}
\item[[1]] P. John, H. Sachs and H. Zernitz. Problem 5: Domino covers in
 square chessboards. Zastosowania Matematyki (Applicationes Mathematicae),
 XIX, 3-4 (1987), 636-638.
\item[[2]] P.W. Kasteleyn, Graph Theory and Crystal Physics, in: F. Harary
 (ed.), Graph Theory and Theoretical Physics, Academic Press, London -
 New York 1967.
\end{itemize}
\end{document}